\input AHTOH-E.STY
\hfuzz4.1pt

\UDC{
512.543.72 
%
}

\MSC{
20F70,   
20F65,   
57M07    
}

\title{%
Yet another Freiheitssatz:
Mating finite groups with locally indicable ones
}
\author{
Anton A. Klyachko
\quad
Mikhail A. Mikheenko
}
\address{
\myAddressWC
\quad
mamikheenko@mail.ru
}
\grants{%
\RSF 22-11-00075}

\abstract{
The main result includes as special cases
\-
on the one hand, the Gerstenhaber--Rothaus theorem (1962)
and its generalisation due to
Nitsche and Thom (2022)
\-
and, on the other hand,
the Brodskii--Howie--Short theorem (1980--1984)
generalising Magnus's Freiheitssatz (1930).
}

\s 0.
Introduction

An equation
$w(x,y,\dots)=1$ over a group $G$
(where $w(x,y,\dots)$ is a word in the
alphabet $G\sqcup\{x^{\pm1},y^{\pm1},\dots\}$) is called
\emph{solvable over} $G$ if some overgroup $\~G\supseteq G$ contains elements
$\~x,\~y,\dots$ such that $w(\~x,\~y,\dots)=1$.

The study of this notion has a long history,
see, e.g.,
[GR62],
[Le62],
[Ly80],
[How81],
[B84],
[EH91],
[How91],
[K93],
[KP95],
[FeR96],
[K97],
[K99],
[CG00],
[EdJu00],
[IK00],
[Juh\'a03],
[K06],
[BK12],
[KT17],
[BE18],
[ABA21],
[EH21],
[NT22],
and references therein.
The following result is
well known.

\proclaim Gerstenhaber--Rothaus theorem
{\rm(for a single equation) [GR62] (see also [LS80])}.
If an equation $g_1x^{\epsilon_1}\dots g_nx^{\epsilon_n}=1$
\(with~unknown $x$\)
over a finite group $G\ni g_1,\dots,g_n$ is \emph{nonsingular},
i.e. $\sum\epsilon_i\ne0$, then it is solvable over~$G$.

It is unknown whether the same holds for any (infinite) groups;
this is
the (strengthened)
Kervaire--Laudenbach conjecture.
Pestov [P08] showed that
this holds true for
hyperlinear groups.%
\fn{%
Definitions, examples, and properties of hyperlinear
(= Connes-embeddable) groups can be found, e.g., in [T18];
we note
only that
the class of hyperlinear groups contains
all finite group and their
free products (possibly, even all group
are hyperlinear
--- this is a well-known
open question).}
A further generalisation
was obtained quite recently.

\proclaim Nitsche--Thom theorem {\rm[NT22] (for a single equation)}.
An equation
$v(x,y,\dots)=1$
over any finite
\(and even any hyperlinear\) group $G$
is solvable over $G$
if the content of $v(x,y,\dots)$ is
nontrivial.

Here, the \emph{content} of a word $v$
from the free product $G*F(x,y,\dots)$
of a group $G$ and the free group~$F(x,y,\dots)$
is the image of $v$ under the natural homomorphism
$\epsilon\:G*F(x,y,\dots)\to F(x,y,\dots)$
(whose kernel is the normal closure~$\nc G$
of $G$).

The following Freiheitssatz
for
\emph{locally indicable}
groups
(i.e. groups whose nontrivial finitely generated subgroups
admit epimorphisms onto $\Z$)
is a result of quite another sort.

\proclaim Brodskii--Howie--Short theorem
{\rm[B80], [B84], [How82], [Sh81]}.
The natural mappings~%
$
C\to(C*D)/\nc w\ot D
$
are injective
if the groups $C$ and $D$ are locally indicable,
and the word $w\in C*D$ is not conjugate to an element
of~$C\cup D$.

\noindent
This generalisation
of Magnus's Freiheitssatz [Mag30]
can easily be reformulated in the
language of equations:
\disp{\sl
any non-exotic equation over a locally indicable group
is solvable over it
}%
(where an equation $w(x,y,\dots)=1$ over a group $G$
is called \emph{exotic} if the word
$w(x,y,\dots)\in G*F(x,y,\dots)$ is conjugate
to an element of $G$).
And conversely: the Nitsche--Thom theorem can be reformulated as
a Freiheitssatz:
\disp{\sl
if $K$ is a free product of finite groups
\(or, more generally, $K$ is hyperlinear\), $F$ is a
free group, and
the image of a word $w\in K*F$ under
the natural homomorphism $K*F\to F$
is nontrivial, then the natural mapping $K\to(K*F)/\nc w$ is
injective.
}%
However, no relations between the Brodskii--Howie--Short and
Gerstenhaber--Rothaus theorems were
known so far.
The purpose of this paper is to fill this lacuna,
i.e. ``to mate" finite and locally indicable groups.
The following fact includes all the results mentioned above as special
cases.

\proclaim Freiheitssatz.
If groups
$C$ and $D$ are locally indicable,
$K$ is a GR$^*$-group,
and
the image of a word~$w\in C*D*K$
under the natural homomorphism
$C*D*K\to C*D$ is
not conjugate to an element of~$C\cup D$,
then the natural mappings
$C*K\to(C*D*K)/\nc w\ot D*K$
are
injective.

\noindent
The definition of
a GR$^*$-group can be found in the next section;
examples of GR$^*$-groups are all hyperlinear groups,
in particular, all free products of finite
groups. Therefore,
this theorem contains the mentioned-above results:
\-
if
$K=\1$,
we obtain the Brodskii--Howie--Short theorem;
\-
if
$C=F(x,y,\dots)$ and $D=F(x_1,y_1,\dots)$ are free groups,
and $w=v(xx_1,yy_1,\dots)$,
\newline
where $v(x',y',\dots)\in(K*F(x',y',\dots)\setminus\nc K$,
then
we obtain the Nitsche--Thom theorem (for a single equation), generalising
the Gerstenhaber--Rothaus theorem (for a single equation).

\enditem
See Section~2 for the complete statement of our main result
and Section 4 for the proof (while Section 3 contains the key lemma).

The full statements of the Gerstenhaber--Rothaus
and Nitsche--Thom
theorems dealing with systems of equations
can be found in Sections 1 and 3.
We do not try to generalise these full versions
in this paper
(because the Brodskii--Howie--Short theorem is essentially
one-relator), but we use these results
(moreover, our approach is heavily based on ideas from [NT22]).
The last section contains a proof of the
Nitsche--Thom theorem (which was essentially obtained in
[NT22], but there are some nuances, see the last section).

The full version of the Brodskii--Howie--Short theorem
is also stronger than the above
statement (and we use it), namely,
one can add
the following words to
the statement
above: \emph{``\dots\ Moreover the group $(C*D)/\nc w$ is
locally indicable too, if the word $w$ is not a proper power"}
(and the following phrase to the
``equational" version of this theorem: \emph{``\dots, and a solution can
be found in a locally indicable overgroup of the given group"}).

The authors thank the Theoretical Physics and Mathematics
Advancement Foundation ``BASIS".

\s 1.
GR- and GR$^*$-group

A system of equations over a group is called \emph{nonsingular}
if the integer rows composed of the exponent-sums of unknowns
in equations are linearly independent.
For example, the system
(with unknowns $x,y,z,t$ over a group~$G\ni a,b,c,d$):
$$
\left\{
\eqalign{
axbycyz^5dz^{-2}&=1
\cr
[xt,dz]^{\the\year}dx^4cy^5bz^6&=1
\cr
ax^7y^8dz^k&=1
}
\right.
\qqbox{has the exponent-sum matrix}
\pmatrix{
1&2&3&0
\cr
4&5&6&0
\cr
7&8&k&0
\cr
},
$$
i.e. this system is singular if and only if $k=9$.

\proclaim Gerstenhaber--Rothaus theorem {\rm[GR62]}.
Any nonsingular system of equations over a finite group is solvable over
this group.

We suggest to
call a group $G$ a \emph{GR-group}, if
any nonsingular system of equations over $G$ is solvable over $G$.
Pestov [P08] noticed that all hyperlinear group are GR-groups.

It is easy to see that
\disp{\sl
the class of GR-groups is a quasivariety,
}%
i.e. this class consists of all groups
satisfying some (possibly, infinite)
system of \emph{quasi-identities},
i.e. (finite) formulae of the form
$$
(\forall x,y,\dots)
\bigl(u_1(x,y,\dots)=1\;\&\dots\;\&\;u_k(x,y,\dots)=1
\imp v(x,y,\dots)=1\bigr),
$$
where $u_i(x,y,\dots)$ and $v(x,y,\dots)$ are words in the alphabet
$\{x^{\pm1},y^{\pm1},\dots\}$.

A simple way to prove this is to apply
the following characterisation of quasivarieties [Mal70]:

\disp{\sl\narrower
a nonempty class of groups is a
quasivariety if and only if
it is closed under passage to subgroups and
reduced \(= filtered\) products.
}%
For the class of GR-groups the both conditions hold obviously.

\medskip

A disadvantage of the quasivariety of GR-groups is
that it is unclear whether this class is closed under
free products. We call a group $G$ a \emph{GR$^*$-group} if it
satisfies the following equivalent conditions:

\item{1)}
the free product $G*\Z$ is a GR-group;

\item{2)}
the free product $G*G$ is a GR-group;

\item{2$'$)}
the free product of
any family of groups isomorphic to $G$
is a GR-group;

\item{3)}
$G$ embeds into a GR-group $\~G\supseteq G$
satisfying no nontrivial
\emph{mixed identity with constants from $G$}, i.e.,
for any word $v(x,y,\dots)\in G*F(x,y,\dots)\setminus\1$
(where $F(x,y,\dots)$ is a free group),
there exist $\~x_v,\~y_v,\dots\in\~G$ such that
$v(\~x,\~y,\dots)\ne1$.

\Proposition.
These conditions are indeed equivalent.

\Proof
First, note that, if $|G|\le2$, then all the conditions are
equivalent, because they are true in this case
(the class of residually finite groups is closed with respect to
free products [Gr57]
and is contained in the class of GR-groups by the Gerstenhaber--Rothaus
theorem, therefore, if $|G|=2$, then we can take,
e.g., the free product of $G$ and the free group of countable rank
as $\~G$ in 3)).

If  $|G|>2$, then 1) and 2) are equivalent, because
$G*\Z$ and $G*G$ embed into each other%
\fn{%
Henceforth, we use the following well-known fact,
which we leave to the reader as an easy exercise:
{\sl
subgroups of a non-dihedral free product
are described up to isomorphism as follows:
these are all groups
allowed by the
cardinality
restriction
and
the Kurosh subgroup theorem.
}}
(and any subgroup of a GR-group is a GR-group).

Condition 2$'$) is equivalent to 1) and 2) by similar reasons:
\-
the free product of any family of groups isomorphic to
$G$ is residually a finite free product $G*\dots*G$
(because any element of an arbitrary free product lies in a 
finite subproduct, and each subproduct is a retract of the whole product);
\-
GR is a residual property (because the class of GR groups is a 
quasi-variety);
\-
a finite free product $G*\dots*G$ (with at least two factors)
and $G*\Z$ embed into each other
by the same ``converse of the Kurosh theorem"*$^)$.

\enditem  
Moreover, Condition 1) implies Condition 3),
because $G*\Z$ contains the free product
$\~G=G*F_\infty$ of $G$ and the free group of infinite (countable)
rank, which, obviously, has no mixed identities with constants
from $G$.

It remains to prove the implication $3)\imp1)$.
The Cartesian product
$H=\bigtimes\limits_{v\in(\~G*\gp t_\infty)\setminus\1}\~G_v$
of copies $\~G_v$ of $\~G$ is a GR-group
(because the Cartesian product of any family of GR-groups is, obviously,
a GR-group). The group $G*\gp t_\infty$
embeds into $H$ as follows: $G$ embeds diagonally,
and the element $t$ is mapped to the element
whose $v$th coordinate is $\~t_v\in\~G$
(i.e. $v(\~t_v)\ne 1$).
Clearly, this is an embedding. It remains to recall that a subgroup of a
GR-group is a GR-group too. This completes the proof.

\Question 1.
Are the classes of GR- and GR$^*$-groups closed with respect to direct and
free products?

\noindent
Simple considerations show that
\-
{\sl the class
of GR-groups is, obviously, closed with respect to direct products};
\-
{\sl if
the class of GR-groups is closed under free products,
then GR-groups and GR$^*$-groups are the same}
(this follows immediately from Condition 2)
of the
definition of GR$^*$-groups;
more generally, any GR-group
decomposable nontrivially into a
free product
is a GR$^*$-group);
\-
{\sl if the class of GR$^*$-groups is closed under direct
products \(or, more generally, if there exists
a GR$^*$-group containing
any two given GR$^*$-groups
as
subgroups\), then the class of GR$^*$-groups
is closed under free products}.

\noindent
This means that the answers to Question 1 can be only the following:
\-
either GR=GR$^*$,
and this class is closed with respect to both operations,
\-
or the class of GR-groups is closed under direct products, but
not closed under free products,
and the class of GR$^*$-groups is
\itemitem{--}
either closed under
neither operations,
\itemitem{--}
or closed with respect to both operations,
\itemitem{--}
or closed under free products, but not closed
under direct products.

\noindent
Certainly,
none of the presently known facts contradicts the
equalities
GR=GR$^*$=$\{$all groups$\}$ (Howie's conjecture).
The class of hyperlinear groups is closed under free products
[BDJ08], hence, all hyperlinear groups are not only GR-
but also GR$^*$-groups.

\s 2. Main result

Suppose that $A\nin G$ and $F$ is a free group.
We call the natural epimorphism $\epsilon_A\:G*F\to(G/A)*F$
the \emph{$G/A$-content};
the $G/G$-content is just the \emph{content} $\epsilon$
(defined in
[KT17] and [NT22]).

\proclaim{Main theorem}.
Suppose that a group $G$ contains a normal subgroup $A$,
which is GR-group,
and the quotient group~$G/A$ is locally indicable.
Then an equation $w(x,y,\dots)=1$
is solvable over $G$ if
the $G/A$-content of $w(x,y,\dots)$ is
not conjugate to an element of $G/A$
in $(G/A)*F(x,y,\dots)$.

\noindent
This theorem implies immediately the Freiheitssatz
stated in the introduction.
Indeed,
consider
the group~$G=C*D*K$,
its normal subgroup
$A=\nc K$
(which
is GR-group,
because it is
isomorphic to the free product of a family of
groups isomorphic to the
GR$^*$-group $K$),
and the embedding $\phi\:G\to G*\gp t_\infty$, where
$C\ni c\mapsto c^t\in G*\gp t_\infty$
and the other free factors
are mapped identically.
The equation
$\phi(w)=1$ (with one variable $t$) is solvable over $G$ by the
main theorem, i.e., the composition
$G\too^\phi G*\gp t_\infty\to G*\gp t_\infty/\nc{\phi(w)}$
is injective on
$D*K$
and contains $w$ in its kernel, as required.
The injectivity of the natural mapping of~$C*K$ can be proven
similarly.

\s 3.
Key lemma

A proof of the following result (obtained essentially in [NT22]) can be
found in the last section.

\proclaim Nitsche--Thom theorem.
A system of equations
$
\{w_1(x,y,\dots)=1,\;w_2(x,y,\dots)=1,\dots\}
$
over a GR-group $G$ is solvable over $G$
if the standard complex of the presentation
$
\pres<x,y,\dots|
\epsilon(w_1)=1,\;\epsilon(w_2)=1,\dots>
$
\(where $\epsilon$ stands for the content\)
admits a covering
with trivial second homologies \(over $\Z$\).

To apply this theorem, we need the following lemma.

\Lemma.
Suppose that a locally indicable group $L$ acts freely on a set
$X$, and $L\times F(X)\too^\o F(X)$ is the natural extension
of this
action to the free group with basis $X$. Then, for each
word $v\in F(X)\setminus\1$, some covering of the standard complex $K$
of the
presentation $H=\pres<X|L\o v>$ has trivial second homologies.

\Proof
Suppose that $v=u^k$ and the word $u$ is not a proper power in $F(X)$.
Let us prove that
the
required covering is the covering $p\:\~K\to K$ corresponding
to the normal closure $\nc{L\o u}\nin H=\pi_1(K)$
of the
orbit of $u$.
In an explicit form,
the complex $\~K$ is the Cayley graph of the group
$A=\pres<X|L\o u>$
with ``$k$ times wrapped" 2-disks
glued to
each cycle with label~$l\o u$, where
$l\in L$ (so, the words $l\o v=l\o u^k$
are written on the
boundaries of the disks).

The group
$A=\pres<X|L\o u>$ embeds into $\^A=Y*L/\nc{\^u}$,
where $Y$ is a set of representatives of orbits of the action
of $L$ on $X$,
and the word $\^u$ is obtained from $u$ by replacing each letter
$x$ with $lyl^{-1}$, where $l\in L$ and $y\in Y$ are the (only) elements
such
that $x=l\o y$. The embedding $A\to\^A$ is obvious:
$x=l\o y\mapsto lyl^{-1}$. The group $\^A$ is locally indicable [B84] and
$\^A=L\semitimes A$.
This group $\^A$ acts
on the complex $\~K$:
the
action of $A$ is the standard action of a group on its
Cayley graph, and $L$ acts by conjugation on vertices: $l\o a=lal^{-1}$
(then vertices joined by an edge labeled $x$ are mapped to vertices
joined by an edge labelled $l\o x$).
This action of $\^A$ is
\-
transitive on two-cells and on vertices
\-
and free on (oriented) edges.

\enditem
The first is quite obvious; and the second can be explained easily:
the stabiliser of an edge must stabilise its origin;
the action is vertex-transitive; therefore, it suffices to show
that the stabiliser of each edge $e$ outgoing from the identity is
trivial; so, the equality $(l,a)\o e=e$ implies that
$a=1$ (otherwise, we would obtain an edge with another origin); now,
the
equality $l\o e=e$ means that $l=1$, because the action of $L$ on $X$
is free.

Now, the situation is simple.
Since the action on two-cells is transitive,
the nontriviality of the second homologies
implies an equality $r\=u=0$, where
\-
$\=u$ is the sum of edges of the boundary of a 2-cell
(which does not vanish, because the universal covering
of a group with one relator $u=1$,
which is not a proper power, is contractible,
i.e. one-relator torsion-free groups are aspherical,
see [LS80]),
\-
and $r$ is a nonzero element of the group ring $\Z[\^A]$ of $\^A$.

\enditem
The freeness of the action on edges leads immediately to zero divisors
in the group ring $\Z[\^A]$
of the locally indicable group $\^A$, which is a contradiction [Hig40].

\s 4.
Proof of the main theorem

Since the word $\epsilon_A(w)\in(G/A)*F(x,y,\dots)$ is not conjugate
to an element of $G/A$, the equation
$\epsilon_A(w)=1$ has a solution in a locally indicable
group $B$
containing $G/A$ (by the Brodskii--Howie--Short theorem).
Therefore,
we have the embeddings
of the $G$ into the unrestricted wreath products:
$
G\subseteq A\Wr(G/A)\subseteq A\Wr B=A^B\leftsemitimes B
$
(where the first embedding is the Kaloujnine--Krasner theorem
[KK51], see also [KaM82]).
Therefore, replacing
$G$ with $A\Wr B$,
and
$A$ with its Cartesian power $A^B\nin A\Wr B$
(which is also a GR-group),
we can assume that the equation $\epsilon_A(w)=1$
has a solution 
$\^x,\^y,\dots\in G/A$; 
hence, making an obvious change of variables 
($x\mapsto x\^x,\ y\mapsto y\^y,\dots $),
we obtain that $\epsilon_A(w)$ is contained in the normal
closure~$\nc F$ of~$F=F(x,y,\dots)$ in $(G/A)*F$.
Thus, we assume
that $w=w(x,y,\dots)$ can be rewritten as a word $\^w$ in the alphabet
$\{x^{\pm b},y^{\pm b},\dots\;|\;b\in B=G/A\}\sqcup A$
(because we assume already that $G=A\leftsemitimes B$).

We arrive to the situation of
the Lemma.
For $X=\{bxb^{-1},byb^{-1},\dots\;|\;b\in B\}$
and $v=\epsilon(\^w)\in F(X)$,
where $\epsilon\:F(X)*A\to F(X)$ is the
natural retraction, and $L=B$,
the Lemma says that a covering of the standard complex
of the presentation~$H=\pres<X|B\o\epsilon(\^w)>$ has trivial second
homologies. By the Nitsche--Thom theorem, this means the solvability
of the system of
equations $\{b\o w=1\;|\;b\in B\}$ (with unknowns $X=B\o\{x,y,\dots\}$)
over the group $A$
(where the action of $B$ on~$X$ is extended to
$A*F(X)$ naturally: $b\o a\:=bab^{-1}$ for~$a\in A$).
In other words, the natural mapping
$A\to\~A=\bigl(A*F(X)\bigr)/\nc{\{b\o w=1\;|\;b\in B\}}$ is injective.
The group $B$ acts on~$\~A$, and this action
extends the action of $B$ on $A$ by conjugations.
Therefore, the natural mapping
$G=A\leftsemitimes B\to\~G=\~A\leftsemitimes B$ is injective.
Moreover, $w(x,y,\dots)=1$ in $\~G$, which means the solvability
of the equation.

\s 5.
Subgroup presentations and the Nitsche--Thom theorem

In this section, we prove the Nitsche--Thom theorem as stated in
Section~3. This result was essentially proven in~[NT22], but
unfortunately, a weaker version was explicitly stated there (see [NT22],
Theorem 1.3 and Remark 2.2).  Our approach somewhat differs from the proof
that can be extracted from [NT22].

The following fact is well known,
see, e.g., [ZVC88],
Theorem 2.2.1, or [LS80], Proposition II.4.1.

\proclaim Schreier's theorem on subgroup presentations {\rm[Sch27]}.
Suppose that $U$ is a subgroup of a group $G=\pres<X|R>$, and
$T\subseteq F(X)$ is a Schreier system of the right coset representatives
of $\pi^{-1}(U)$ in the free group $F(X)$, where $\pi\:F(X)\to G=F/\nc R$
is the natural epimorphism. Then the subgroup $U$ is generated by the
images of the words~$y_{t,x}=tx(\={tx})^{-1}\in F(X)$ and has a
presentation, where generators are all nontrivial
words $\{y_{t,x}\;|\;t\in T, x\in X\}$, and the defining relators are all
words $\{trt^{-1}\;|\;t\in T, x\in X\}$ rewritten as words in nontrivial
generators~$y_{t,x}$ \(considered as letters\).

\noindent
Here, the bar means taking a representative:
$\=v\in T$ (where $v\in F(X)$) is the unique word from~$T$
such that $\pi^{-1}(U)\=v=\pi^{-1}(U)v$;
and $T$ is \emph{Schreier} in the usual sense:
any prefix of any word from~$T$ lies in $T$.

\noindent
The geometric interpretation of this
presentation, which we call \emph{Schreier},
is also well known:
\-
the group $G$ is the fundamental group of the standard complex
$K$ corresponding to the presentation $\pres<X|R>$
(i.e. $K$ has a single vertex,
edges correspond to generators from $X$,
and two-dimensional cells correspond to relators from $R$);
\-
the subgroup $U$ is the fundamental group of the complex $\~K$
(with a base vertex) covering the complex~$K$;
\-
vertices of $\~K$ correspond to right cosets of
$U$ in $G$ (and
of $\pi^{-1}(U)$
in $F(X)$);
\-
a Schreier system of representatives corresponds to a maximal
subtree in the 1-skeleton of $\~K$;
\-
generators $y_{t,x}$ correspond to edges of $\~K$
(or, more precisely, each edge $e\in\~K$
corresponds to a path starting at the base vertex
going through the maximal subtree to the start-point of $e$,
passing $e$, and returning to the base vertex via the maximal
subtree); so, nontrivial generators $y_{t,x}$ correspond to
edges not belonging to the maximal subtree;
\-
finally, relators of $U$ correspond to 2-cells of $\~K$.

\enditem
This geometric interpretation shows
that
\disp{\sl
the standard complex of the presentation $G=\pres<X|R>$
admits a covering with trivial second homologies if and only if
the relators of the Schreier presentation
of some subgroup $U\subseteq G$ form a nonsingular
system \(in the sense of Section 1\).
}%
Thus, we can restate the Nitsche--Thom theorem
in group-theoretical terms
(with no topology or homologies).

\proclaim Nitsche--Thom theorem {\rm(a pure group-theoretical form)}.
A system of equations $W=1$
\(possibly, infinite and, possibly, with an infinite set
of unknowns $X$\) over a GR-group $G$
is solvable over $G$ if
there exists a subgroup of~$\pres<X|\epsilon(W)>$,
the relators of whose Schreier
presentation form a nonsingular system.
{\rm (Here, $\epsilon$ stands for the content, see Introduction.)}

\Proof
Suppose that $A\subseteq\pres<X|\epsilon(W)>$ is a subgroup
with a nonsingular Schreier presentation.
Consider a presentation~$G=\pres<Z|V>$
of the group $G$,
the corresponding presentation of the group
$\bigl(G*F(X)\bigr)/\nc W=\pres<Z\sqcup X|V\sqcup W>=H$,
the natural epimorphism
$\theta\:H\to H/\nc G=\pres<X|\epsilon(W)>$,
and the subgroup $\theta^{-1}(A)\subseteq H$.
It is easy to see that the Schreier presentation
for the subgroup~$\theta^{-1}(A)\subseteq H$
transforms into the Schreier presentation
for the subgroup~$A\subseteq\pres<X|\epsilon(W)>$ by deleting
all generators corresponding to conjugate of generators from $Z$
(in geometrical terms,
each vertex of the
covering complex corresponding
to the subgroup~$\theta^{-1}(A)\subseteq H$
is contained in a subcomplex
isomorphic to the
standard complex of the presentation~$G=\pres<Z|V>$;
and, contracting each such subcomplex to a single point,
we obtain
the covering complex corresponding to the subgroup
$A\subseteq\pres<X|\epsilon(W)>$).

\noindent
Thus, the relators of
the group
$\theta^{-1}(A)$ form
a nonsingular system of equations
over the free product of
$\Bigl|\pres<X|\epsilon(W)>:A\Bigr|$ copies of $G$.
Since $G$ is a GR-group,
the natural mapping
$G\to\theta^{-1}(A)\subseteq H$ is injective%
\fn{%
but we cannot assert that the entire
free product of $\Bigl|\pres<X|\epsilon(W)>:A\Bigr|$ copies of
$G$ embeds into $\theta^{-1}(A)$; it embeds if $G$ is not only
GR- but also a GR$^*$-group.}
(to show this, we can take the quotient group of $\theta^{-1}(A)$ by
the normal closure of the product of all copies of $G$, except
one). This implies immediately the injectivity of the natural
homomorphism $G\to H$, as required.

\baselineskip 10.2pt
\References

[ABA21]
M. F. Anwar, M. Bibi, M. S. Akram,
On solvability of certain equations of arbitrary length
over torsion-free groups,
Glasgow Mathematical Journal, 63:3 (2021), 651-659.
\arXiv 1903.06503

[BK12]
D. V. Baranov, A. A. Klyachko,
Economical adjunction of square roots to groups,
Siberian Math. Journal, 53:2 (2012), 201-206.
\arXiv 1101.3019

[B80]
S. D. Brodskii,
Equations over groups and group with a single defining relation,
Russian Math. Surveys, 35:4 (1980), 165-165.


[B84]
S. D. Brodskii,
equations over groups and group with one relator,
{Siberian Math. Journal}, 25:2 (1984), 235-251.

[BDJ08]
N. Brown, K. Dykema, K. Jung,
Free entropy dimension in amalgamated free products,
Proc. London Math. Soc., 97:2 (2008), 339-367.
\arXiv math/0609080

[BE18]
M. Bibi, M. Edjvet,
Solving equations of length seven over torsion-free groups,
Journal of Group Theory, 21:1 (2018), 147-164.

[CG00]
A. Clifford, R. Z. Goldstein,
Equations with torsion-free coefficients,
{Proc. Edinburgh Math. Soc.}, {43:2} (2000), 295-307.

[EH91]
M. Edjvet, J. Howie,
The solution of length four equations over groups,
{Trans. Amer. Math.  Soc.}, {326:1} (1991), 345-369.

[EH21]
M. Edjvet, J. Howie,
On singular equations over torsion-free groups,
International Journal of Algebra and Computation, 31:3 (2021), 551-580.
\arXiv:2001.07634

[EdJu00]
M. Edjvet, A. Juh{\accent 19 a}sz,
Equations of length 4 and one-relator products,
{Math. Proc. Cambridge Phil. Soc.}, 129:2 (2000), 217-230.

[FeR96]
R. Fenn, C. Rourke,
Klyachko's methods and the solution of
equations over torsion-free groups,
{L'Enseignment
Math{\accent 19 e}matique}, 42 (1996), 49-74.

[GR62]
M. Gerstenhaber, O.S. Rothaus,
The solution of sets of equations in groups,
{Proc. Nat. Acad. Sci. USA}, {48:9} (1962), 1531-1533.

[Gr57]
K. W. Gruenberg,
Residual properties of infinite soluble groups,
Proc. London Math. Soc., s3-7:1 (1957), 29-62.

[Hig40]
G. Higman,
The units of group-rings,
Proc. London Math. Soc., s2-46:1 (1940), 231-248.

[How81]
J. Howie,
On pairs of 2-complexes and systems of equations over groups,
{J. Reine Angew Math.}, 1981:324 (1981), 165-174.

[How91]
J. Howie,
The quotient of a free product of groups by a single high-powered relator.
III: The word problem,
{Proc. Lond. Math. Soc.}, {62}:3 (1991), 590-606.

[IK00]
S. V. Ivanov, A. A. Klyachko,
Solving equations of length at most six over torsion-free groups,
Journal of Group Theory, 3:3 (2000), 329-337.

[Juh{\accent 19 a}03]
Juh{\accent 19 a}sz A.
On the solvability of a class
of equations over groups,
{Math. Proc. Cambridge Phil. Soc.}, 135:2 (2003), 211-217.

[KaM82]
M. I. Kargapolov, Ju. I. Merzljakov,
Fundamentals of the theory of groups.
Graduate Texts in Mathematics, 62, Springer, 1979.

[K93]
A. A. Klyachko,
A funny property of sphere and equations over groups,
Communications in Algebra, 21:7 (1993), 2555-2575.

[K97]
A. A. Klyachko,
Asphericity tests,
International Journal of Algebra and Computation, 7:4 (1997), 415-431.

[K99]
A. A. Klyachko,
Equations over groups, quasivarieties,
and a residual property of a free group,
Journal of Group Theory, 2:3 (1999), 319-327.

[K06]
A. A. Klyachko,
How to generalize known results on equations over groups,
Math. Notes, 79:3 (2006), 377-386.
\arXiv math.GR/0406382

[KP95]
A. A. Klyachko, M. I. Prishchepov,
The descent method for equations over groups,
Moscow Univ. Math. Bull. 50:4 (1995), 56-58.

[KT17]
A. A. Klyachko, A. B. Thom,
New topological methods to solve equations over groups,
Algebraic and Geometric Topology, 17:1 (2017), 331-353.
\arXiv 1509.01376

[KK51]
M. Krasner, L. Kaloujnine,
Produit complet des groupes de permutations et le
probl\`eme d'extension de groupes. III,
Acta Sci. Math., 14 (1951), 69-82.

[Le62]
F. Levin,
Solutions of equations over groups,
{Bull. Amer. Math. Soc.}, {68:6} (1962), 603-604.

[Ly80]
R. C. Lyndon,
Equations in groups,
{Bol. Soc. Bras. Math.}, {11}:1 (1980), 79-102.

[LS80]
R. Lyndon, P. Schupp,
Combinatorial group theory.
Springer, 1977.

[Mag30]
W. Magnus,
{\accent "7F U}ber diskontinuierliche Gruppen mit einer
definierenden Relation (Der Freiheitssatz),
{J. Reine Angew Math.}, {163} (1930) 141-165.

[Mal70]
A. I. Mal'cev,
Algebraic systems.
Springer, 1973.

[NT22]
M. Nitsche, A. Thom,
Universal solvability of group equations,
Journal of Group Theory, 25:1 (2022), 1-10.
\arXiv 1811.07737

[P08]
V. G. Pestov,
Hyperlinear and sofic groups: A brief guide,
Bull. Symb. Log., 14:4 (2008), 449-480.
\arXiv:0804.3968

[Sch27]
O. Schreier,
Die Untergruppen der freien Gruppen,
Abhandlungen aus dem Mathematischen Seminar der Universit\"at Hamburg,
5:1 (1927), 161-183.

[Sh81]
H. Short,
Topological methods in group theory: the adjunction problem.
Ph.D. Thesis,
University of Warwick, 1981.

[T18]
A. Thom,
Finitary approximations of groups and their applications,
in: Proceedings of the International Congress of Mathematicians ---
Rio de Janeiro 2018. Vol. III.
Invited lectures, World Scientific, Hackensack (2018), 1779-1799.
\arXiv 1712.01052

[ZVC88]
H. Zieschang, E. Vogt, H.-D. Coldewey,
Surfaces and planar discontinuous groups.
Springer (1980).

\end